\documentclass[a4paper]{amsart}
\usepackage{amssymb,latexsym,amsmath}
\usepackage{amsmath}
\usepackage{enumerate}
\theoremstyle{plain}

\theoremstyle{definition}

\theoremstyle{definition}

\numberwithin{equation}{section}
\newtheorem{Theorem}{\quad Theorem}[section] 

\newtheorem{Definition}[Theorem]{\quad Definition} 

\newtheorem{Corollary}[Theorem]{\quad Corollary} 

\newtheorem{Lemma}[Theorem]{\quad Lemma} 

\newtheorem{Proposition}[Theorem]{\quad Proposition} 

\newtheorem{Remark}[Theorem]{\quad Remark}

\begin{document}
\title[ Some Properties of Accretive Operators in Linear 2-Normed Spaces]{Some Properties of Accretive Operators in Linear 2-Normed Spaces }

\author[P.K. Harikrishnan]{P.K. Harikrishnan}

\address{Department of   Mathematics,  Manipal Institute of Technology, Manipal University, Manipal-576104, Karnataka, India.}

\email{pk.harikrishnan@manipal.edu, pkharikrishnans@gmail.com }

\author[K.T. Ravindran]{K.T. Ravindran}
\address{Department of Mathematics, Payyanur College, Kannur University, Kannur, India}

\email {drktravindran@gmail.com}
\begin{abstract}
In this paper we discuss some properties of resolvents of an accretive operator in linear 2-normed spaces, focusing on the concept of contraction mapping and the unique fixed point of contraction mappings in linear 2- normed spaces. Also, we establish the existence of solution of strong accretive operator equation in linear 2-normed spaces.
\end{abstract}

\subjclass[2000]{41A65, 41A15}

\keywords{ Linear 2-normed spaces, sequentially closed, accretive operators, m-accretive operator, fixed point, sequentially continuous }

\maketitle

\section{Introduction} The concept of  2- metric spaces, linear 2- normed spaces and 2-inner product spaces , introduced by S Gahler in 1963 , paved the way for a number of authors like, A White, Y J Cho, R Freese, C R Diminnie, to do work on possible applications of Metric geometry, Functional Analysis and Topology as a new tool. A systematic presentation of the recent results related to the Geometry of linear 2-normed spaces as well as an extensive list of the related references can be found in the book [1]. In [4] S Gahler introduced the following definition of linear 2-normed spaces. The study of accretive operators in various spaces has been done by the authors T Kato[6], Shih-sen Chang, Yeol Je Cho, Shin Min Kang [5].. 

\section{Preliminaries} 

\begin{Definition}[4] Let X be a real linear space of dimension greater than 1 and $\|.,.\|$  be a real valued function on $X\times X$ satisfying the properties,
\newline \textbf{A1:}$\|x,y\|$ = 0 iff x and y are linearly dependent
\newline \textbf{A2:}$\|x,y\|$ = $\|y,x\|$
\newline \textbf{A3:}$\|\alpha x,y\|$ = $|\alpha|$$\|y,x\|$
\newline \textbf{A4:}$\|x+y,z\|$ $\leq$ $\|x,z\|+ \|y,z\|$

then the function $\|.,.\|$ is called a 2-norm on X. The pair (X,$\|.,.\|$)  is called a linear 2- normed space.
\end{Definition} 

Some of the basic properties of 2-norms, they are non-negative and $\|x,y+\alpha x\|$ = $\|x,y\|$ $\forall$ x,y $\in$ X and $\alpha \in$ R.

The most standard example for a linear 2 - normed space is X = R$^{2}$ equipped with the following 2-norm,
\begin{center}
$\|x_{1},x_{2}\|$ = abs $\det$$\left(
                            \begin{array}{cc}
                              x_{11} & x_{12}\\
                              x_{21} & x_{22} \\
                            \end{array}
                          \right)$ where $x_{i}=( x_{i1}, x_{i2})$for i=1,2
\end{center}

Every linear 2-normed space is a locally convex TVS. In fact, for a fixed b$\in$ X, P$_{b}$(x)= $\|x,b\|$ for x $\in$ X is a semi norm  and the family $\{P_{b};b \in X\}$ of semi norms generates a locally convex topology on X.

\begin{Definition}[1] Let $(X,\|.,.\|)$  be a linear 2- normed space, E be a subset of X then the closure of E is $\overline{E}=\{x\in X ;$ there is a sequence $x_{n}$ of E such that $x_{n}\rightarrow x $.
We say, E is sequentially closed if $E=\overline{E}$.
\end{Definition} 

\begin{Definition}[1] Let $(X,\|.,.\|)$  be a linear 2- normed space, then the mapping $T: X\rightarrow X$ is said to be sequentially continuous if $x_{n}\rightarrow x$ implies $Tx_{n}\rightarrow Tx$.
\end{Definition}

\begin{Definition}[3]  An operator $T:D(T)\subset X\rightarrow X$ is said to be accretive if for every $z\in D(T)$
\begin{center}
    $\|x-y,z\|\leq\|(x-y)+\lambda (Tx-Ty),z\|$ for all $x,y \in D(T)$ and $\lambda > 0$.
\end{center}
\end{Definition} 

\begin{Definition}[3] An operator $T:D(T)\subset X\rightarrow X$ is said to be m-accretive if $R(I+\lambda T) = X$ for $\lambda > 0$.
\end{Definition} 

\begin{Remark}[4]
Every closed subset of a 2-Banach space is complete.
\end{Remark}

\section{Main Results}Consider an accretive operator  $T:D(T)\subset X\rightarrow X$  then for  $\lambda > 0$ ,$(I+\lambda T)^{-1}$ exists.We have,$T:D(T)\rightarrow R(T)$  then  $(I+\lambda T)$  is onto. Let$ x_{1},x_{2}\in D(T) $ with $(I+\lambda T)x_{1} $ = $(I+\lambda T)x_{2}$ implies  $ (x_{1}-x_{2})+\lambda (Tx_{1}-Tx_{2})=0 $ implies $ \|(x_{1}-x_{2})+\lambda (Tx_{1}-Tx_{2}),z\|\ = 0 $ for every $ z\in X $ . Since T is accretive , we have  $\Vert x_{1}-x_{2},z\Vert \leq 0$ for every  $ z\in X $ implies $\Vert x_{1}-x_{2},z\Vert = 0$ for every  $ z\in X $ implies $ x_{1}-x_{2} =0 $  implies $ x_{1}= x_{2} $ . So $(I+\lambda T)$  is one-one.\\ Hence for  $\lambda > 0$ , $(I+\lambda T)^{-1}$exist.

\begin{Definition} Let $(X,\|.,.\|)$  be a linear 2- normed space then an operator T on X is said to be non expansive if for each x,y $ \in $D(T), $ \|Tx-Ty,z\| \leq \|x-y,z\| $ for every z$ \in $X. T is said to be expansive if for each x,y $ \in $D(T), $ \|Tx-Ty,z\| > \|x-y,z\| $ for every z$ \in $X.
\end{Definition} 

\begin{Proposition} .Let $(X,\|.,.\|)$  be a linear 2- normed space and $T:D(T)\subset X\rightarrow X$ be an operator. If $(I+ \lambda T)$ is expansive for all $\lambda > 0$ then T is accretive.
\end{Proposition} 
Proof: The proof is immediate from the definition 3.1.

\begin{Remark}
If T is an accretive operator on a 2-normed space X then $(I+\lambda T)^{-1}$ is non expansive.
\end{Remark}

\begin{Definition} Let $(X,\|.,.\|)$  be a linear 2- normed space E be a non empty subset of X and $ e \in E $ then E is said to be e-bounded if there exists some $ M > 0 $ such that $ \|x,e\|\leq M $ for all $ x \in E $. If for all $ e \in E $, E is e-bounded then E is called a bounded set.
\end{Definition}

Let $(X,\|.,.\|)$be a linear 2- normed space and T be an m-accretive operator on X. For n=1,2,3.., Define the resolvent of T as, $ J_{n}(x)=(I+n^{-1}T)^{-1}(x) $ and the Yosida approximation $ T_{n}(x)=n(I-J_{n})(x) $ for all $ x \in X $ and $ \lambda > 0 $.

\begin{Proposition} .Let $(X,\|.,.\|)$  be a linear 2- normed space and T be a m-operator then

    (i)$\|J_{n}x-J_{n}y,z\|\leq\|(x-y),z\|$ 
    
    (ii)$\|T_{n}x-T_{n}y,z\|\leq 2n\|(x-y),z\|$ for all $x,y \in X$
\end{Proposition} 
 
Proof:
 
 (i) Since T is accretive in X,we have $(I+\alpha T)^{-1}$ is non expansive 
 
 implies for each x,y $ \in $D(T), $ \|(I+\alpha T)^{-1}x-(I+\alpha T)^{-1}y,z\| \leq \|x-y,z\| $ for every z$ \in $X 
 
 Take $ \alpha = 1/n $ then $ \alpha > 0 $. So $\|J_{n}x-J_{n}y,z\|\leq\|(x-y),z\|$.
 
 (ii) We have for each x,y $ \in $D(T), 
\begin{center}
$\|T_{n}x-T_{n}y,z\|=\|n(I-J_{n})x-(I-J_{n})y,z\|$ 

$~~~~~~~~~~~~~~~~~~~$ $=\|n(x-y)+n(J_{n}x-J_{n}y),z\|$ 

$~~~~~~~~~~~~~~~~~~~$ $\leq n\Vert x-y,z \Vert + n\Vert J_{n}x-J_{n}y,z\Vert$

$~~~~~~~~~~~~~~~~~$ $\leq n\Vert x-y,z \Vert + n\Vert x-y,z \Vert$

$~~~~~~~~~~~~~~~$ $= 2n\Vert x-y,z \Vert$

\end{center}
for all $ z \in X $

\begin{Proposition} .Let $(X,\|.,.\|)$  be a linear 2- normed space and E be a non empty bounded subset of X and $T:E\rightarrow E$ be an accretive operator then there exists some $ M>0 $ and for $ x \in E $ , $ \Vert T_{n}x,z\Vert \leq M $ for all $ z \in E $.
\end{Proposition}
 
 Proof: Let $ x \in E $ then 
\begin{center}
 $ T_{n}x=n(x-J_{n}x)=n(J_{n}(J^{-1}_{n}x)-J_{n}x)=n(J_{n}(I+n^{-1}T)x-J_{n}x) $
\end{center}

So, by Proposition 3.5 (i);
\begin{center}
 $ \Vert T_{n}x,z\Vert =\Vert n(J_{n}(I+n^{-1}T)x-J_{n}x),z \Vert \leq n\Vert (I+n^{-1}T)x-x,z \Vert  =\Vert Tx,z\Vert \leq M $
\end{center}
for all $ z \in E $
 
 \begin{Proposition} .Let $(X,\|.,.\|)$  be a linear 2- normed space and E be a non empty bounded subset of X and $T:E\rightarrow E$ be an accretive operator . If $ x \in \overline{E} $ then $ J_{n}x \longrightarrow x $.
\end{Proposition}
 
 Proof: Let $ x \in E $ then$ \Vert x-J_{n}x,z\Vert = n^{-1}\Vert T_{n}x,z \Vert \leq n^{-1}M $ for some $ M>0 $ and for all $ z \in E $.
 
 ie; $ \Vert x-J_{n}x,z\Vert \longrightarrow 0$ as $ n\longrightarrow \infty $ implies $ J_{n}x \longrightarrow x $

 If  x $\in \overline {E} $ then there exists a sequence $ \lbrace x_{m}\rbrace \in E $ such that $ x_{m}\longrightarrow x $.
 Since $ x_{m} \in E $ we have $ J_{n}(x_{m}) \longrightarrow x_{m} $.
 
 We have, for all $ z \in E $

 $ \Vert J_{n}x-x,z\Vert =\Vert J_{n}x-J_{n}(x_{m})+J_{n}(x_{m})-x_{m} +x_{m}-x,z \Vert $ 
 
 $~~~~~~~~~~~~~~~~$ $\leq \Vert J_{n}x-J_{n}(x_{m}),z \Vert + \Vert J_{n}(x_{m})-x_{m},z \Vert +\Vert x_{m}-x,z\Vert $
 
 $~~~~~~~~~~~~~~~~$ $\leq \Vert x-x_{m},z \Vert + \Vert J_{n}(x_{m})-x_{m},z \Vert +\Vert x_{m}-x,z\Vert \longrightarrow 0$ as $ n\longrightarrow\infty $
 
 Hence $ J_{n}x \longrightarrow x $.

\begin{Definition} Let $(X,\|.,.\|)$  be a linear 2- normed space then the mapping $ T:X\longrightarrow X $ is said to be a contraction if there exists some $ k\in(0,1) $ such that $ \Vert Tx-Ty,z \Vert \leq k \Vert x-y,z \Vert $ for all x,y,z $ \in $ X.
\end{Definition}

\begin{Lemma} Let $(X,\|.,.\|)$  be a linear 2- normed space then every contraction $ T:X\rightarrow X $ is sequentially continuous.
\end{Lemma} 

Proof: Since T is a contraction, there exists some $ k\in(0,1) $ such that $ \Vert Tx-Ty,z \Vert \leq k \Vert x-y,z \Vert $ for all x,y,z $ \in $ X.

Let $ \lbrace x_{n}\rbrace $ be a sequence in X such that $ x_{n}\longrightarrow x $.

Then, $ \Vert Tx_{n}-Tx,z \Vert \leq k \Vert x_{n}-x,z \Vert \longrightarrow 0$ as $ n\longrightarrow \infty $ implies 
$ Tx_{n}\longrightarrow Tx $. So T is sequentially continuous.

Next lemma proves the analogues of Banach fixed point theorem in Linear 2-normed spaces.

\begin{Lemma} Let $(X,\|.,.\|)$  be a linear 2- normed space and E be a non empty closed and bounded subset of X. Let  $ T:E\rightarrow E $ be a contraction then T has a unique fixed point on X.
\end{Lemma} 

Proof:  Since T is a contraction, there exists some $ k\in(0,1) $ such that 
\begin{center}
$ \Vert Tx-Ty,z \Vert \leq k \Vert x-y,z \Vert $ for all x,y,z $ \in $ X.
\end{center}
We have, $ \Vert T^{2}x-T^{2}y,z \Vert = \Vert T(Tx)-T(Ty),z \Vert \leq k \Vert Tx-Ty,z \Vert \leq k^{2} \Vert x-y,z \Vert $  for all $ z \in E $. Similarly $ \Vert T^{n}x-T^{n}y,z \Vert  \leq k^{n} \Vert x-y,z \Vert $ for all $ z \in E $.

Let $ x_{0} \in E $ then construct a sequence depending on $ x_{0} $.

Let $ x_{1}=Tx_{0}, x_{2}=Tx_{1}, x_{3}=Tx_{2} ,..., x_{n}=Tx_{n-1} $ then $ x_{1}=Tx_{0}, x_{2}=T^{2}x_{0}, x_{3}=T^{3}x_{0} ,..., x_{n}=T^{n}x_{0} $.

We show that $ \lbrace x_{n} \rbrace $ is a cauchy sequence in E. Let $ m,n > 0$ with $m>n$. Take $m=n+p$ then for any $z \in E$,

$\Vert x_{n}-x_{m},z\Vert = \Vert x_{n} - x_{n+p},z\Vert$

$ ~~~~~~~~~~~~~~~~~ $ $ = \Vert (x_{n} - x_{n+1})+ (x_{n+1}- x_{n+2})+ x_{n+2}-...+ (x_{n+p-1}- x_{n+p}),z\Vert$

$ ~~~~~~~~~~~~~~~~~ $ $ \leq \Vert x_{n} - x_{n+1},z\Vert +\Vert x_{n+1}- x_{n+2},z\Vert +...+ \Vert x_{n+p-1}- x_{n+p},z\Vert $

$ ~~~~~~~~~~~~~~~~~~~~~~ $ $ =  \Vert T^{n}x_{0} -T^{n}x_{1},z\Vert +\Vert T^{n+1} x_{0}- T^{n+1}x_{1},z\Vert +...+ \Vert T^{n+p-1}x_{0}- T^{n+p-1}x_{1},z\Vert $

$ ~~~~~~~~~~~~~~~~~ $ $ \leq k^{n}\Vert x_{0} -x_{1},z\Vert +k^{n+1}\Vert x_{0} -x_{1},z\Vert +...+k^{n+p-1}\Vert x_{0} -x_{1},z\Vert $

$ ~~~~~~~~~~~~~~~~~ $ $ \leq k^{n}\Vert x_{0} -x_{1},z\Vert (1+k+k^{2}+...) = (\frac{k^{n}}{1-k})\Vert x_{0} -x_{1},z\Vert$

1e; $\Vert x_{n}-x_{m},z\Vert \leq (\frac{k^{n}}{1-k})\Vert x_{0} -x_{1},z\Vert$ for all $ z \in E $

Since E is bounded, there exists $M>0$ such that $ \Vert x_{0} -x_{1},z\Vert \leq M $ for all $ z \in E $

So, $\Vert x_{n}-x_{m},z\Vert \leq (\frac{k^{n}M}{1-k})$ for all $ z \in E $ implies $\Vert x_{n}-x_{m},z\Vert \longrightarrow 0$ as $n\longrightarrow\infty$, since $k \in (0,1)$ imples $ \lbrace x_{n} \rbrace $ is a cauchy sequence in E, So $ \lbrace x_{n} \rbrace $  converges to some x in E.

Since T is sequentially continuous, $ Tx = \lim Tx_{n} = \lim x_{n+1}= x $ as $n\longrightarrow\infty$. Therefore T has a fixed point in E. 

Now, we have to prove that such a fixed point is unique.

Let $ y \in E $ with $ y \neq x $ and $ Ty \neq y $.

Assume that $ Ty=y $ then $ \Vert x-y,z \Vert = \Vert Tx-Ty,z\Vert \leq k\Vert x-y,z \Vert $ implies $ k\geq 1 $, a contradiction to $ k \in (0,1) $. So, $ Ty \neq y $. Hence the fixed point of T is unique in E.

\begin{Definition} Let $(X,\|.,.\|)$  be a linear 2- normed space then the mapping $T:X\rightarrow X$ is said to be strong accretive if for every $z\in X$
\begin{center}
    $\|(\lambda -k)(x-y),z\|\leq\|(\lambda -1)(x-y)+(Tx-Ty),z\|$ 
\end{center}
for all $x,y \in X$ and $\lambda > k,$ and $ k \in (0,1) $.
\end{Definition}

\begin{Theorem}Let $(X,\|.,.\|)$  be a linear 2- normed space and E be a non empty closed and bounded subset of X. Let  $ T:E\rightarrow X $ be a strong accretive mapping. If $ (I+T)(E) \supset E $ then the equation $ Tx=\theta $ has a solution in E, where $ \theta $ is the zero vector in X.
\end{Theorem}

Proof: Since T is strong accretive implies 
\begin{center}
$\|(\lambda -k)(x-y),z\|\leq\|(\lambda -1)(x-y)+(Tx-Ty),z\|$  
\end{center}

for all $x,y \in X$ and $\lambda > k,$ and $ k \in (0,1)$

Put $ \lambda=2 $ then $\|(2 -k)(x-y),z\|\leq\|(x-y)+(Tx-Ty),z\|$  

implies  $\vert 2-k \vert\|(x-y),z\|\leq\|(I+T)(x) - (I+T)(y)),z\|$ 

implies $\|(x-y),z\|\leq \frac{\|(I+T)(x) - (I+T)(y),z\|}{\vert 2-k \vert}$ 

Since $ (I+T)(E) \supset E $ we get,

 $ x,y \in D[(I+T)^{-1}] $ implies $ x= (I+T)^{-1}(u)$ and $ y= (I+T)^{-1}(v)$ for some $u,v \in E$

So, $\|(I+T)^{-1}(u)-(I+T)^{-1}(v),z\|\leq \frac{\| u - v,z\|}{\vert 2-k \vert}$

ie; $ (I+T)^{-1} : E\rightarrow E $ is a contraction mapping. By the Lemma 3.10 $ (I+T)^{-1} $ has a fixed point $ x_{0} $ in E.
ie; $ x_{0}=(I+T)^{-1}(x_{0}) $ implies $ x_{0}=(I+T)(x_{0}) $ implies $ Tx_{0} = \theta$.

Hence,the equation $ Tx=\theta $ has a solution in E.

\begin{Corollary}Let $(X,\|.,.\|)$  be a linear 2- normed space and E be a non empty closed and bounded subset of X. Let  $ T:E\rightarrow E $ be a strong accretive mapping. If $ (I+T)(E) = E $ then $R(T)=E$.
\end{Corollary}

Proof:

For any $ p \in E $. Take $ T_{0}=T-p $

Since T is strong accretive we have,

$\|(\lambda -k)(x-y),z\|\leq\|(\lambda -1)(x-y)+(Tx-Ty),z\|$  for all $x,y \in X$ and $\lambda > k,$ and $ k \in (0,1) $

We have, $ T_{0} $ is also strong accretive.

Therefore, by theorem 3.12, there exists $ x_{0}\in E $ such that $ T_{0}(x_{0})= \theta $ implies $ Tx_{0}-p=\theta $ implies $ Tx_{0}=p $ for every $ p \in E $. Therefore, $R(T)=E$.

\end{document}